\documentclass[a4paper,12pt]{article}
\usepackage{amssymb}

\newcommand{\R}{\mathbb{R}}

\newcommand{\N}{\mathbb{N}}

\newcommand{\beq}{\begin{equation} }
\newcommand{\eqq}{\end{equation} }
\newcommand{\cuad}{{\sqcap\kern-.68em\sqcup}}

\newtheorem{teo}{Theorem}[section]

\newtheorem{proposition}{Proposition}[section]

\newtheorem{lemma}{Lemma}[section]
\newtheorem{corollary}{Corollary}[section]
\newtheorem{remark}{Remark}[section]
\newcommand{\bremark}{\begin{remark} \em}
\newcommand{\eremark}{\end{remark} }

\def\beeq{\begin{equation}}
\def\eeq{\end{equation}}
\newcommand{\begeqaet}{\begin{eqnarray*}}
\newcommand{\eneqaet}{\end{eqnarray*}}

\begin{document}
\begin{center}{\bf  \Large  Singular solutions of fractional elliptic equations with absorption }\medskip
\medskip

\bigskip

{\bf Huyuan Chen $^1$}

 Departamento de Ingenier\'{\i}a  Matem\'atica
\\ Universidad de Chile,
  Santiago, Chile\\[1mm]

{\bf Laurent V\'{e}ron $^2$}

Laboratoire de Math\'{e}matique et Physique Th\'{e}orique
\\ CNRS UMR 7350
\\  Universit\'{e} Fran\c{c}ois Rabelais, Tours, France\\[1mm]

 {\sl  ($^1$hchen@dim.uchile.cl, $^2$Laurent.Veron@lmpt.univ-tours.fr)}

\bigskip
\bigskip
\begin{abstract}
The aim of this paper is to study the singular solutions to
fractional elliptic equations with absorption
$$
\left\{ \arraycolsep=1pt
\begin{array}{lll}
 (-\Delta)^\alpha  u+|u|^{p-1}u=0,\quad & \rm{in}\quad\Omega\setminus\{0\},\\[2mm]
u=0,\quad & \rm{in}\quad \R^N\setminus\Omega,\\[2mm]
\lim_{x\to 0}u(x)=+\infty,
\end{array}
\right.
$$
where $p>0$, $\Omega$ is an open, bounded and smooth  domain of
$\R^N\ (N\ge2)$ with $0\in\Omega$.

We analyze the existence, nonexistence, uniqueness and asymptotic
behavior of the solutions.

\end{abstract}

\bigskip


\end{center}

\setcounter{equation}{0}
\section{Introduction}
In the present paper, we  are concerned with the singular solutions  of fractional elliptic problems of the form
\begin{equation}\label{eq1.1}
\left\{ \arraycolsep=1pt
\begin{array}{lll}
 (-\Delta)^\alpha  u+|u|^{p-1}u=0,\quad & \rm{in}\quad\Omega\setminus\{0\},\\[2mm]
u=0,\quad & \rm{in}\quad \R^N\setminus\Omega,\\[2mm]
\lim_{x\to 0}u(x)=+\infty,
\end{array}
\right.
\end{equation}
where  $\Omega$ is an open, bounded and smooth  domain
of $\R^N\ (N\ge2)$ with $0\in\Omega$, $p>0$ and $(-\Delta)^\alpha $  with $\alpha\in(0,1)$ is the fractional Laplacian
defined as
\begin{equation}\label{eq 1}
(-\Delta)^\alpha u(x)=P.V.\int_{\R^N}\frac{u(x)-u(y)}{|x-y|^{N+2\alpha}}dy,
\end{equation}
here $P.V.$ stands for the principle value integral, that for notational simplicity we omit in what follows.

During the last years, singular solutions of nonlinear elliptic
equations have been studied by many authors. We just mention the
earlier work by V\'{e}ron \cite{LV1,LV2}, Gmira-V\'{e}ron \cite{GV},
Brezis-Lions \cite{BL}, Bandle-Marcus \cite{BM1,BM2}, Baras-Pierre
\cite{BP}, Chen-Matano-V\'{e}ron \cite{CMV}, without any attempt to
review the references here. The first result of unconditional
removability of isolated sets for semilinear elliptic equations with
absorption term is due to Brezis-V\'{e}ron \cite{BV}. They
considered the classical equation
\begin{equation}\label{eqjt1.1}
-\Delta   u+g(u)=0\quad  \rm{in}\quad\Omega\setminus\{0\},
\end{equation}
where $\Omega$ is an open subset of $\R^N\ (N\ge3)$ containing $0$
and $g$ is a continuous function satisfying some extra hypothesis,
then there exists a solution for equation (\ref{eqjt1.1}) in the
whole $\Omega$. Later on, this result was extended in \cite{LV2}
using the method which is developed by Baras-Pierre \cite{BP}. In
the meantime,  V\'{e}ron \cite{LV1} has done much work for equation
(\ref{eqjt1.1}) with $g(u)=|u|^{p-1}u$ and $1<p<\frac{N}{N-2}$ if $
N\ge3$ $(p>1 \ if\ N=2)$, he described the behaviour of solution for
equation (\ref{eqjt1.1}) near the isolated singularity.

Recently, great attention has been devoted to investigate nonlinear
equations involving fractional Laplacian. Caffarelli-Silvestre
\cite{CS} gave a new formulation of the fractional Laplacian through
Dirichlet-Neumann maps. Later, they studied the regularity results
for fractional problems in \cite{CS1,CS2}. The existence of solution
for equation with fractional Laplacian was proved by Cabr\'{e}-Tan
\cite{CT}, Felmer-Quaas \cite{FQ}, Servadei-Valdinoci \cite{SV}.
Moreover, Li \cite{L}, Chen-Li-Ou \cite{CLO1,CLO2} and Felmer-Wang
\cite{FW} studied symmetry results and monotonicity of positive
solutions for fractional equations.  Chen-Felmer-Quaas \cite{CFQ}
analyzed the existence and asymptotic behavior of  large solution to
fractional equation with absorption by advanced   method of super
and sub solutions.

The purpose of this paper is to study singular solutions for
fractional equations  (\ref{eq1.1}) with absorption, including the
existence and the asymptotic behavior of singular solutions near
$0$. It is well-known that the singular near 0 of functions
$|x|^\tau$ only could be considered with $\tau\in(-N,0)$ for working
by fractional laplacian,  which is a nonlocal operator. In the
following, we  state  main result.

\begin{teo}\label{teo 1}
Suppose  that  $\Omega$ is an open,  bounded and smooth domain of
$\R^N\ (N\ge2)$ with $0\in\Omega$, $\alpha\in(0,1)$.\\[2mm]
$(i)$ If
\begin{equation}\label{1.4}
1+\frac{2\alpha}N<p<\frac N{N-2\alpha},
\end{equation}
then  problem (\ref{eq1.1}) admits a positive solution $u$ such that
for some $C>0$,
\begin{equation}\label{1.3}
 \lim_{x\to 0}u(x)|x|^{\frac{2\alpha}{p-1}}=C.
 \end{equation}
 Moreover, that solution $u$ is unique in the sense of
\begin{equation}\label{1.2}
 0<\liminf_{x\to 0}u(x)|x|^{\frac{2\alpha}{p-1}}\le \limsup_{x\to
0}u(x)|x|^{\frac{2\alpha}{p-1}}<+\infty.
 \end{equation}
$(ii)$ If
\begin{equation}\label{1.5}
0<p<\frac N{N-2\alpha},
\end{equation}
then for any $t>0$,  problem (\ref{eq1.1}) admits a positive
solution $u$ such that
\begin{equation}\label{1.7}
\lim_{x\to0}u(x)|x|^{N-2\alpha}=t.
 \end{equation}

\noindent$(iii)$ If $p>0$, then problem (\ref{eq1.1}) doesn't admit
any solution $u$ such that
\begin{equation}\label{1.70}
 0<\liminf_{x\to 0}u(x)|x|^{-\tau}\le \limsup_{x\to
0}u(x)|x|^{-\tau}<+\infty,
 \end{equation}
for any $\tau\in(-N,0)\setminus\{2\alpha-N,-\frac{2\alpha}{p-1}\}$.

\end{teo}

Theorem \ref{teo 1} part $(i)$ presents the existence, uniqueness in
the sense of (\ref{1.2}) and the asymptotic behavior with power
$-\frac{2\alpha}{p-1}$ of singular solution to (\ref{eq1.1}), part
$(ii)$ shows the existence and the asymptotic behavior with power
$-N+2\alpha$ of singular solution to (\ref{eq1.1}) and part $(iii)$
gives the nonexistence of singular solution to (\ref{eq1.1}) in the
sense (\ref{1.70}). In the next, we give some remarks to show more
information for singular solution to (\ref{eq1.1}).

\begin{remark}\label{re 0}
Under the hypothesis of Theorem \ref{teo 1} part $(i)$, the solution
$u$, which satisfies (\ref{1.3}), has estimate
\begin{equation}
|u(x)-C_1|x|^{-\frac{2\alpha}{p-1}}|<C_2,\quad
x\in\Omega\setminus\{0\}
 \end{equation}
where  $C_1>0$ will be given in (\ref{3.4.1}) and $C_2>0$.
\end{remark}

\begin{remark}\label{re 1}
Under the hypothesis of Theorem \ref{teo 1} part $(ii)$, if
\begin{equation}\label{1.8}
\frac{2\alpha}{N-2\alpha}<p<\frac N{N-2\alpha},
\end{equation}
then for any $t>0$,  problem (\ref{eq1.1}) admits a positive
solution $u$ such that, for any $0<|x|<d_0$, we have
\begin{equation}\label{1.9}
\frac{|x|^{\tau_1}}C\le t|x|^{2\alpha-N}-u(x)\le C|x|^{\tau_1},
 \end{equation}
 where $C>0$, $\tau_1=2\alpha-(N-2\alpha)p<0$ and
$d_0=\frac{1}3 \min\{dist(0,\partial\Omega),1\}$.
\end{remark}

\begin{remark}\label{re 2}
Under the hypothesis of Theorem \ref{teo 1} part $(iii)$,
 if $p\ge \frac N{N-2\alpha}$, then problem
(\ref{eq1.1}) doesn't admit any solution $u$ such that
\begin{equation}\label{1.1}
 0<\liminf_{x\to 0}u(x)|x|^{-\tau}\le \limsup_{x\to
0}u(x)|x|^{-\tau}<+\infty,
 \end{equation}
for any $\tau\in(-N,0).$
\end{remark}

The rest of the paper is organized as follows. In Section \S 2, we
introduce Preliminaries for existence and some estimates which is
used for constructing super and sub solutions of (\ref{eq1.1}). In
Section \S3, we prove the existence of the solutions of
(\ref{eq1.1}). The uniqueness is addressed in Section \S4. In the
section \S 5, it is devoted to non-existence.


\setcounter{equation}{0}
\section{Preliminaries}

We remind here  some basic knowledge about $(-\Delta)^\alpha$ with
$\alpha\in(0,1)$, see for instance \cite{CFQ}.
\begin{lemma}\label{cr 2.0}
 Assume that $x_0$ achieves the maximum of $u$ in $\R^N$, then
 \begin{equation}\label{3.2.1}
 (-\Delta)^{\alpha} u(x_0)\ge 0.
 \end{equation}
Moreover, if $x_0$ achieves the maximum of $u$ in $\R^N$, then
\begin{equation}\label{3.2.1}
 (-\Delta)^{\alpha} u(x_0)\le 0,
 \end{equation}
 holds if and only if
$$u(x)=u(x_0)\ \  \rm{a.e.\ in}\ \ \R^N.$$
\end{lemma}

\begin{lemma}\label{th 2.1}
Assume that $0\in\Omega$ and $p>0$. Moreover, we suppose that there
are super-solution $\bar U$ and sub-solution $\underline{U}$
of (\ref{eq1.1}) such that
$$\bar U\geq \underline{U}\ \ \rm{in}\ \Omega\setminus\{0\},
\quad \liminf_{x\to 0}\underline U(x)=+\infty,\quad \bar
U=\underline{U}=0\ \ \rm{in}\ \Omega^c.$$ Then there exists at least
one positive solution $u$ of (\ref{eq1.1}) such that
$$\underline{U}\leq u\leq \bar U\ \ \rm{in}\
\Omega\setminus\{0\}.$$
\end{lemma}
{\bf Proof.}  The process of the proof is the same as Theorem 2.6 in
\cite{CFQ}.
\hfill$\Box$\\

In order to construct super and sub solutions for problem
(\ref{eq1.1}), we will use some appropriate truncated functions. To
describe our following analysis, we give some notations. By
$0\in\Omega$,  it is able to assume that $\delta\in(0,d_0)$ is such
that $d(\cdot)=dist(\cdot,\partial\Omega)$ is  $C^2$ in
$A_\delta:=\{x\in \Omega\ |\ d(x)<\delta\}$ and $d(x)\le |x|$ in
$A_\delta$, where $d_0=\frac{1}3 \min\{dist(0,\partial\Omega),1\}$.
Let $B_r:=B_r(0)\setminus\{0\}$ for any $r>0$, we have
$dist(A_\delta, B_{d_0})>0$. Moreover, we define
\begin{equation}\label{2.1}
V_\tau(x):=\left\{ \arraycolsep=1pt
\begin{array}{lll}
 |x|^{\tau},\ & x\in B_{d_0},\\[2mm]
 d(x)^{2},\ & x\in A_\delta,\\[2mm]
 l(x),\ \ \ \ &
x\in \Omega\setminus  (A_\delta\cup  B_{d_0}(0)),\\[2mm]
0,\ &x\in\Omega^c,
\end{array}
\right.
\end{equation}
 where  $\tau$ is a parameter in $(-N,0)$ and the function $l$  is  positive  such
 that $l(x)\le |x|^\tau$ in $\Omega\setminus  (A_\delta\cup  B_{d_0}(0))$ and
$V_\tau$ is $C^2$ in $\R^N\setminus\{0\}$.

It will be convenient for next auxiliary lemmas
to define the following function
 \begin{equation}\label{2.2}
C(\tau):= \int_{\R^N}\frac{|z-e_1|^\tau-1}{|z|^{N+2\alpha}}dz
\end{equation}
where $e_1=(1,0,\cdots,0)\in\R^N$.
It is well known from
\cite{FQ} that
\begin{equation}\label{2.4}
C(\tau)\left\{ \arraycolsep=1pt
\begin{array}{lll}
 >0,\quad &\rm{if}\quad \tau\in(-N,-N+2\alpha),\\[2mm]
 =0,\ &\rm{if}\quad   \tau=-N+2\alpha,\\[2mm]
 <0,\ \ \ \ &\rm{if}\quad
\tau\in(-N+2\alpha,0).
\end{array}
\right.
\end{equation}

\begin{lemma}\label{prop 2.1}
Assume that  $\Omega$ is an open, bounded, smooth domain with
$0\in\Omega$ and $\tau\in(-N,0)$. Then there exists $c>0$ such htat
\begin{equation}\label{2.30}
-c<(-\Delta)^\alpha V_\tau(x)+C(\tau)|x|^{\tau-2\alpha}\le0, \ \ \
\forall\ x\in B_{d_0/2},
\end{equation}
where $C(\cdot)$ is defined in (\ref{2.2}).
\end{lemma}
{\bf Proof.}  For any given $x\in B_{d_0/2}$, we have
\begin{eqnarray*} -(-\Delta)^\alpha V_\tau(x)
&=&\int_{\R^N}\frac{V_\tau(z)-V_\tau(x)}{|z-x|^{N+2\alpha}}dz
=\int_{\R^N}\frac{V_\tau(z)-|x|^\tau}{|z-x|^{N+2\alpha}}dz
\\&=&\int_{\R^N}\frac{|z|^\tau-|x|^\tau}{|z-x|^{N+2\alpha}}dz+\int_{\R^N\setminus
B_{d_0}}\frac{V_\tau(z)-|z|^\tau}{|z-x|^{N+2\alpha}}dz
\\&=:&I_1(x)+I_2(x).
\end{eqnarray*}
We look at each of these integrals separately. On one side, by direct computation, we have
$$I_1(x)=\int_{\R^N}\frac{|z+x|^\tau-|x|^\tau}{|z|^{N+2\alpha}}dz=C(\tau)|x|^{\tau-2\alpha}.$$
On the other side, for $z\in \R^N\setminus B_{d_0}$ and $x\in B_{d_0/2}$,
we have $|z-x|\geq\frac{|z|}{2}$ and $|V_\tau(z)-|z|^\tau|\leq c|z|^\tau$ for some $c>0$. Then there exists $C>0$ such that
\begin{eqnarray*}
I_2(x)&=&\int_{\R^N\setminus
B_{d_0}}\frac{V_\tau(z)-|z|^\tau}{|z-x|^{N+2\alpha}}dz
\\&\ge &-C\int_{\R^N\setminus
B_{d_0}}|z|^{\tau-N-2\alpha}dz\\& \ge& -Cd_0^{\tau-2\alpha}.
\end{eqnarray*}
On the other hand by $V_\tau(z)\le |z|^\tau$, we have
$$I_2(x)=\int_{\R^N\setminus
B_{d_0}}\frac{V_\tau(z)-|z|^\tau}{|z-x|^{N+2\alpha}}dz<0.$$  Hence,
we obtain (\ref{2.30}). The proof is compete. \hfill$\Box$

As a consequence, we have the following corollary
\begin{corollary}\label{coro 2.1}
Let $\Omega$ be an open, bounded, smooth domain containing
$0$.\\
$(i)$\ If  $$\tau\in(-N,-N+2\alpha),$$ then  there exists
$\delta_1\in(0,d_0)$ and $C>1$ such that
$$
\frac1C|x|^{\tau-2\alpha }\leq-(-\Delta)^{\alpha}V_\tau(x)\leq
C|x|^{\tau-2\alpha },\ \ \forall\ x\in B_{\delta_1}.
$$
 $(ii)$
If $$ \tau\in(-N+2\alpha,0),$$ then  there exists
$\delta_1\in(0,d_0)$ and $C>1$ such that
$$
\frac1C|x|^{\tau-2\alpha }\leq(-\Delta)^{\alpha}V_\tau(x)\leq
C|x|^{\tau-2\alpha },\ \ \forall\ x\in B_{\delta_1}.
$$ $(iii)$\ If
$$\tau=-N+2\alpha,$$
 then
there exists  $C>1$ such that
$$|(-\Delta)^{\alpha}V_\tau(x)|\leq C, \ \ \forall\ x\in\Omega\setminus\{0\}.$$

\end{corollary}
{\bf Proof.} It follows directly Lemma \ref{prop 2.1} and
(\ref{2.4}). \hfill$\Box$


\setcounter{equation}{0}
\section{Existence of Problem (\ref{eq1.1}) }

This section is devoted to use Corollary \ref{coro 2.1} to to
construct suitable sub-solution and super-solution of
(\ref{eq1.1}) to prove the existence.\\[1mm]
{\bf Proof of  Remark \ref{re 0}.} Firstly, we construct
super-solution and sub-solution of (\ref{eq1.1}) under the
hypotheses of Theorem \ref{teo 1} part $(i)$ by adjusting the
parameter $\lambda>0$ in the following functions
\begin{equation}\label{3.4.1}
U_\lambda(x):=C(\tau_p)^{\frac1{p-1}} V_{\tau_p}(x)+\lambda\bar
V(x)\ \ \rm{ and}\ \ W_{\lambda}(x):=C(\tau_p)^{\frac1{p-1}}
V_{\tau_p}(x)-\lambda\bar V(x),
\end{equation}
where $V_{\tau_p}$ is defined in (\ref{2.1}) with
$\tau_p=-\frac{2\alpha}{p-1}\in (-N,-N+2\alpha)$,  $C(\tau_p)>0$ is
defined in (\ref{2.2}) and $\bar V$ is the solution of
\begin{equation} \label{3.1}
\left\{ \arraycolsep=1pt
\begin{array}{lll}
 (-\Delta)^{\alpha} \bar V(x)=1,\ \ \ \ &
x\in\Omega,\\[2mm]
\bar V(x)=0,\ &x\in\Omega^c.
\end{array}
\right.
\end{equation}
By Lemma \ref{cr 2.0}, we have that $\bar V>0$ in $\Omega$.  \\[1mm]
\textbf{1. Super-solution.} By the definition of $U_\lambda$, it has
\begin{eqnarray*} (-\Delta)^\alpha
U_\lambda(x)= C(\tau_p)^{\frac1{p-1}}(-\Delta)^\alpha
V_{\tau_p}(x)+\lambda\quad\rm{in}\
\Omega\setminus\{0\}.\end{eqnarray*}

 By (\ref{2.30}) and $\tau_p p=\tau_p-2\alpha $,   it follows that for all $\lambda\geq
 0$,
\begin{eqnarray*}
(-\Delta)^\alpha U_\lambda(x)+U^p_\lambda(x)\geq -C(\tau_p)^{\frac
p{p-1}} |x|^{\tau_p-2\alpha }+C(\tau_p)^{\frac p{p-1}} |x|^{\tau_p
p}\ge0,\quad x\in B_{\frac{d_0}2}.
\end{eqnarray*}
In above inequality  we used that for any $a,b\ge0$,
$$(a+b)^{p}\geq a^p.$$

Next we consider the domain $ \Omega\setminus B_{\frac{d_0}2}(0)$.
Then, by definition of $V_\tau$, there exists $C_1>0$ such that
$$|(-\Delta)^\alpha V_\tau|\leq C_1\ \ \rm{in}\ \ \Omega\setminus B_{\frac{d_0}2}(0).$$
Then there exists $\bar \lambda>0$ such that for $\lambda\ge \bar
\lambda$, it has
\begin{eqnarray*} (-\Delta)^\alpha
U_\lambda(x)+U_\lambda^p(x)&\geq&\lambda -
C_1C(\tau_p)^{\frac1{p-1}}\ge 0.\end{eqnarray*}

Together with $U_{\bar \lambda}= 0$ in $\Omega^c$, we have that
$U_{\bar\lambda}$ is
a super-solution of  (\ref{eq1.1}).\\[1mm]
\textbf{2. Sub-solution.} We observe that
\begin{eqnarray*}(-\Delta)^\alpha
W_{\lambda}(x)= C(\tau_p)^{\frac1{p-1}}(-\Delta)^\alpha
V_{\tau_p}(x)-\lambda\quad\rm{in}\ \Omega\setminus\{0\}.
\end{eqnarray*}

By (\ref{2.30}),  it follows that for $x\in B_{\frac{d_0}2}$ and
$\lambda\ge0$,
\begin{eqnarray*}(-\Delta)^\alpha
 W_{\lambda}(x)+ |W_{\lambda}|^{p-1}W_{\lambda}(x)\leq
-C(\tau_p)^{\frac p{p-1}} |x|^{\tau_p-2\alpha}+C(\tau_p)^{\frac
p{p-1}} |x|^{\tau_p p}\le0.\end{eqnarray*} In the first inequality
above we used that for any $a,b\ge0$,
$$|a-b|^{p-1}(a-b)\leq a^p.$$

Since $(-\Delta)^\alpha V_\tau+V_\tau^p$ is continuous in $\Omega\setminus\{0\}$,
then there exists $C_2>0$ such that
$$|C(\tau_p)^{\frac 1{p-1}}(-\Delta)^\alpha V_\tau|+ C(\tau_p)^{\frac p{p-1}}V_\tau^p\leq C_2,\ \ x\in
\Omega\setminus B_{\frac{d_0}2}(0).$$
 Then, there exists $\underline{\lambda}>0$ such that for $\lambda\ge \underline{\lambda}$, we have
 \begin{eqnarray*}
(-\Delta)^\alpha
W_{\lambda}(x)+|W_{\lambda}|^{p-1}W_{\lambda}(x)&\leq&
C_2-\lambda\\&\leq&0,\ \ x\in\Omega\setminus B_{\frac{d_0}2}(0).
\end{eqnarray*}
Then $ W_{\underline{\lambda}}$ is a sub-solution of  (\ref{eq1.1}).
Since  $\bar \lambda,\underline{\lambda}>0$ and $\bar V>0$ in
$\Omega$, then
\begin{equation}\label{4.2}
U_{\bar \lambda}> W_{\underline \lambda}\ \ \rm{in}\
\Omega\setminus\{0\}\ \ \rm{and}\ \ U_{\bar \lambda}=
W_{\underline \lambda}=0\ \rm{in}\ \Omega^c.
\end{equation}

Then, by Lemma \ref{th 2.1}, there exists at least one positive
 solution $u$ such that$$ W_{\underline \lambda}\leq u\leq
U_{\bar \lambda}\quad \rm{in}\ \Omega\setminus\{0\}.$$
 The proof is complete.\hfill $\Box$\\

The proof of Theorem \ref{teo 1} part $(i)$ follows the proof of
Remark \ref{re 0}.

\noindent{\bf Proof of Theorem \ref{teo 1} part $(ii)$ with
$0<p\le\frac{2\alpha}{N-2\alpha}$.} Let $\tau_0=2\alpha-N$ and
$\tau_1=\frac{2\alpha-N}2<0$. For $0<p\le\frac{2\alpha}{N-2\alpha}$,
we have that $$0>p\tau_0\ge \tau_1-2\alpha.$$ For any given $t>0$,
we define
$$
U_{\mu}(x):=t V_{\tau_0}(x)+\mu\bar V(x)$$ and
$$ W_{\mu}(x):=t
V_{\tau_0}(x)-\mu V_{\tau_1}(x)-\mu^2\bar V(x),
$$
where $\mu,\lambda>0$, $V_\tau$ is defined in (\ref{2.1}) and $\bar
V$ is the solution of (\ref{3.1}). We construct super-solution and
sub-solution of (\ref{eq1.1}) under the hypotheses of Theorem
\ref{teo 1} part $(ii)$ by adjusting the positive parameters $\mu$.\\[1mm]
\textbf{1. Super-solution.} By the definition of $U_{\mu}$, it has
\begin{eqnarray*} (-\Delta)^\alpha
U_{\mu}(x)= t(-\Delta)^\alpha V_{\tau_0}(x)+\mu,\quad x\in
\Omega\setminus\{0\}.
\end{eqnarray*}
By Corollary \ref{coro 2.1} part $(iii)$, for $x\in B_{d_0}$, it
follows that
\begin{eqnarray*}
(-\Delta)^\alpha U_{\mu}(x)+U^p_{\mu}(x)\geq -Ct +t^p |x|^{\tau_0
p}.\end{eqnarray*} Then there exists $\delta_2\in(0,d_0)$ such that
$$(-\Delta)^\alpha U_{\mu}(x)+U^p_{\mu}(x)\geq 0,\ x\in
B_{\delta_2}.$$

Next we consider the domain $ \Omega\setminus B_{\delta_2}(0)$.
Then, by definition of $U_{\mu,\lambda}$, there exists $C_1>0$ such
that
$$|(-\Delta)^\alpha V_{\tau_0}|\leq C_1\ \ \rm{in}\ \ \Omega\setminus B_{\delta_2}(0).$$
 Then  there exists $\mu_1>1$ such that for
$\mu\ge \mu_1$, it has
\begin{eqnarray*} (-\Delta)^\alpha
U_{\mu}(x)+U_{\mu}^p(x)\geq\mu-tC_1\ge 0,\quad x\in\Omega\setminus
B_{\delta_2}(0).
\end{eqnarray*}
Then  $U_{\mu_1}$ is
a super-solution of  (\ref{eq1.1}).\\[1mm]
\textbf{2. Sub-solution.} We observe that
\begin{eqnarray*}(-\Delta)^\alpha
W_{\mu}(x)= t(-\Delta)^\alpha V_{\tau_0}(x)-\mu(-\Delta)^\alpha
V_{\tau_1}(x)-\mu^2,\quad x\in\Omega\setminus\{0\}.
\end{eqnarray*}
By Corollary \ref{coro 2.1} part $(ii)$ and $(iii)$, for $x\in
B_{\delta_1}$, it follows that
\begin{eqnarray*}(-\Delta)^\alpha
 W_{\mu}(x)+ |W_{\mu}|^{p-1}W_{\mu}(x)&\leq&
Ct-\frac\mu C |x|^{\tau_1-2\alpha }+t^p |x|^{\tau_0
p},\end{eqnarray*} where $C>1$. Here the inequality above we used
that for any $a,b\ge0$,
$$|a-b|^{p-1}(a-b)\leq a^p.$$
Then for $\mu\ge2Ct^p$ and $\tau_1-2\alpha<\tau_0p$, there exists
$\delta_2>0$ such that
$$(-\Delta)^\alpha  W_{\mu}(x)+|W_{\mu}|^{p-1}W_{\mu}(x)\leq 0,\ x\in
B_{\delta_2}.$$

Since $(-\Delta)^\alpha V_\tau+V_\tau^p$ is continuous in
$\Omega\setminus\{0\}$, then there exists $C_2>0$ such that
$$|(-\Delta)^\alpha V_{\tau_0}|+ V_{\tau_0}^p\leq C_2,\ \ x\in
\Omega\setminus B_{\delta_2}(0)$$ and
$$|(-\Delta)^\alpha V_{\tau_1}|\leq C_2,\ \ x\in
\Omega\setminus B_{\delta_2}(0).$$ Then, there exists $\mu_2\ge
2Ct^p$ such that for $\mu\ge \mu_2$, we have
 \begin{eqnarray*}
(-\Delta)^\alpha W_{\mu}(x)+|W_{\mu}|^{p-1}W_{\mu}(x)&\leq& C_2t+\mu
C_2+C_2^pt^p-\mu^2\\&\leq&0,\ \ x\in\Omega\setminus B_{\delta_2}(0).
\end{eqnarray*}
As a consequence, $ W_{\mu_2}$ is a sub-solution of  (\ref{eq1.1}).

Since  $\mu_2,\mu_1>0$ and $\bar V,V_{\tau_0},V_{\tau_1}>0$ in
$\Omega\setminus\{0\}$, then
\begin{equation}\label{4.2}
U_{\mu_1}> W_{\mu_2}\ \ \rm{in}\ \Omega\setminus\{0\}\ \ \rm{and}\ \
U_{\mu_1}= W_{\mu_2}=0\ \rm{in}\ \Omega^c.
\end{equation}

Then by Lemma \ref{th 2.1},  there exists solution  $u$ of
(\ref{eq1.1}) satisfies  (\ref{1.7}). The proof is complete.\hfill
$\Box$\\[2mm]

\noindent{\bf Proof of Remark \ref{re 1}.}
 For any given
$t>0$, we define
$$
U_{\mu,\lambda}(x):=t V_{\tau_0}(x)-\mu
V_{\tau_1}(x)+\lambda\bar V(x)$$
and
$$ W_{\mu}(x):=t
V_{\tau_0}(x)-\mu V_{\tau_1}(x)-\mu^2\bar V(x),
$$
where $\mu,\lambda>0$, $\tau_0=2\alpha-N$, $\tau_1=\tau_0p+2\alpha$,
$V_\tau$ is defined in (\ref{2.1}), and $\bar V$ is the solution of
(\ref{3.1}). By $\frac{2\alpha}{N-2\alpha}<p<\frac{N}{N-2\alpha}$,
we have that $$-N+2\alpha<\tau_1<0.$$  We construct super-solution
and sub-solution of (\ref{eq1.1}) under the hypotheses of Remark
\ref{re 1} by adjusting the positive parameters $\mu$ and
$\lambda$.\\[1mm]
\textbf{1. Super-solution.} By the definition of $U_{\mu,\lambda}$,
it has
\begin{eqnarray*} (-\Delta)^\alpha
U_{\mu,\lambda}(x)= t(-\Delta)^\alpha
V_{\tau_0}(x)-\mu(-\Delta)^\alpha V_{\tau_1}(x)+\lambda,\quad x\in
\Omega\setminus\{0\}.
\end{eqnarray*}
By Corollary \ref{coro 2.1} part $(ii)$ and $(iii)$, for $x\in B_{\delta_1}$,
it follows that
\begin{eqnarray*}
(-\Delta)^\alpha U_{\mu,\lambda}(x)+U^p_{\mu,\lambda}(x)\geq -Ct
-C\mu |x|^{\tau_1-2\alpha }+t^p |x|^{\tau_0 p}.\end{eqnarray*} Then
letting $\mu=t^p/(2C)$ and there exists $\delta_2\in(0,\delta_1)$  such that
$$(-\Delta)^\alpha U_{\mu,\lambda}(x)+U^p_{\mu,\lambda}(x)\geq 0,\ x\in
B_{\delta_2}.$$

Next we consider the domain $ \Omega\setminus B_{\delta_2}(0)$.
Then, by definition of $U_{\mu,\lambda}$, there exists $C_1>0$ such
that
$$|(-\Delta)^\alpha V_\tau|\leq C_1\ \ \rm{in}\ \ \Omega\setminus B_{\delta_2}(0),$$
for $\tau=\tau_0,\tau_1$. Then for $\mu=t^p/(2C) $, there
exists $\lambda_1>1$ such that for $\lambda\ge \lambda_1$, it has
\begin{eqnarray*} (-\Delta)^\alpha
U_{\mu,\lambda}(x)+|U_{\mu,\lambda}|^{p-1}U_{\mu,\lambda}(x)&\geq&\lambda -\mu
C_1-tC_1-\mu^pC_1^p\\&\ge& 0,\quad x\in\Omega\setminus B_{\delta_2}(0).
\end{eqnarray*}

Then  for $ \lambda=\lambda_1>1$ and $\mu=\mu_1=t^p/2$, we have that
$U_{\mu_1,\lambda_1}$ is
a super-solution of  (\ref{eq1.1}).\\[1mm]
\textbf{2. Sub-solution.} We observe that
\begin{eqnarray*}(-\Delta)^\alpha
W_{\mu}(x)= t(-\Delta)^\alpha V_{\tau_0}(x)-\mu(-\Delta)^\alpha
V_{\tau_1}(x)-\mu^2,\quad x\in\Omega\setminus\{0\}.
\end{eqnarray*}
By Corollary \ref{coro 2.1} part $(ii)$ and $(iii)$, for $x\in B_{\delta_1}$, it follows
that
\begin{eqnarray*}(-\Delta)^\alpha
 W_{\mu}(x)+ |W_{\mu}|^{p-1}W_{\mu}(x)&\leq&
Ct-\frac\mu C |x|^{\tau_1-2\alpha }+t^p |x|^{\tau_0
p},\end{eqnarray*} where $C>1$. Here the inequality above we used
that for any $a,b\ge0$,
$$|a-b|^{p-1}(a-b)\leq a^p.$$
Then for $\mu\ge2Ct^p$, there exists $\delta_2>0$ such that
$$(-\Delta)^\alpha  W_{\mu}(x)+|W_{\mu}|^{p-1}W_{\mu}(x)\leq 0,\ x\in
B_{\delta_2}.$$

Since $(-\Delta)^\alpha V_\tau+V_\tau^p$ is continuous in
$\Omega\setminus\{0\}$, then there exists $C_2>0$ such that
$$|(-\Delta)^\alpha V_{\tau_0}|+ V_{\tau_0}^p\leq C_2,\ \ x\in
\Omega\setminus B_{\delta_2}(0)$$ and
$$|(-\Delta)^\alpha V_{\tau_1}|\leq C_2,\ \ x\in
\Omega\setminus B_{\delta_2}(0).$$
Then, there exists $\mu_2\ge 2Ct^p$ such that for $\mu\ge \mu_2$, we have
 \begin{eqnarray*}
(-\Delta)^\alpha W_{\mu}(x)+|W_{\mu}|^{p-1}W_{\mu}(x)&\leq&
C_2t+\mu C_2+C_2^pt^p-\mu^2\\&\leq&0,\ \ x\in\Omega\setminus
B_{\delta_2}(0).
\end{eqnarray*}
As a consequence, $ W_{\mu_2}$ is a sub-solution of  (\ref{eq1.1}).

Since  $\mu_2>\mu_1>0$ and $\bar V,V_{\tau_0},V_{\tau_1}>0$ in $\Omega\setminus\{0\}$, then
\begin{equation}\label{4.2}
U_{\mu_1,\lambda_1}> W_{\mu_2}\ \ \rm{in}\ \Omega\setminus\{0\}\ \ \rm{and}\ \
U_{\mu_1,\lambda_1}= W_{\mu_2}=0\ \rm{in}\ \Omega^c.
\end{equation}

Then by Lemma \ref{th 2.1},  there exists solution  $u$ of
(\ref{eq1.1}) satisfies  (\ref{1.7}). The proof is complete.\hfill
$\Box$

The proof of Theorem \ref{teo 1} part $(ii)$ with
$\frac{2\alpha}{N-2\alpha}<p<\frac{N}{N-2\alpha}$ follows the proof
of Remark \ref{re 1}.



\setcounter{equation}{0}
\section{Proof of the uniqueness}

In this section, we prove the uniqueness in Theorem \ref{teo 1} part
$(i)$ by contradiction.
 Let $u$ and $v$ be two solutions of problem (\ref{eq1.1}) satisfying (\ref{1.2}).
We observe that,  $u$ and $v$ are positive in $\Omega\setminus\{0\}$
and there exists $C_0\ge1$ such that
\begin{equation}\label{4.1.3}
 \frac 1{C_0}\leq v(x)|x|^{-\tau},\ u(x)|x|^{-\tau} \leq C_0,\ \ \forall x\in
  B_{d_0},
\end{equation}
where, we recall, $ B_{d_0}= B_{d_0}(0)\setminus\{0\}$,
$d_0=\frac13dist(0,\partial\Omega)$ and in whole this section,
$\tau=-\frac{2\alpha}{p-1}$ of Theorem \ref{teo 1} part $(i)$. We
denote
\begin{equation}\label{4.1.4}
\mathcal{A}=\{x\in  B_{d_0}|\ u(x)>v(x)\}.
\end{equation}

It is easy to see that $\mathcal{A}$ is open and $\mathcal{A}\subset
\Omega$.

\begin{teo}\label{th 3.1}
Under the hypotheses of  Theorem \ref{teo 1} part $(i)$, we have
$$\mathcal{A}=\O.$$
\end{teo}

To overcome the difficulty caused by the nonlocal character,    we
introduce the following lemmas to prove Theorem \ref{th 3.1}. We
denote  $$g(x)=\left\{ \arraycolsep=1pt
\begin{array}{lll}
 (1-|x|^2)^3,\ \ \ \ &
x\in B_1(0),\\[2mm]
0,\ &x\in B_1^c(0).
\end{array}
\right.$$ Since $g$ is $C^2$  in $\R^N$, then there exists $\bar
C>0$ such that
$$(-\Delta)^\alpha g(x)\leq \bar C,\ \ x\in B_1(0).$$
Then it is obvious to see that
\begin{lemma}\label{lm 2.2}
Let $V=g/\bar C$ in $\R^N$, where $g(x)$ and $\bar C>0$ defined
above, then $$(-\Delta)^\alpha V(x)\leq1$$ and
\begin{equation}\label{2.3} V(0)=\max_{x\in\R^N}V(x).
\end{equation}

 \end{lemma}

\begin{lemma}\label{lm 3.1}
Under the hypotheses of Theorem \ref{teo 1}part $(i)$, if
$$\mathcal{A}_{k,M}:=\{x\in  \R^N\setminus\{0\}\ |\  u(x)-k v(x)> M\}\not=\O,$$
for $k>1$ and $M\ge0$. Then,
\begin{equation}\label{4.1.5}
0\in\partial \mathcal{A}_{k,M}.
\end{equation}
\end{lemma}
{\bf Proof.}  If (\ref{4.1.5}) is not true, there exist $\bar r>0$
such that
$$ \mathcal{A}_{k,M}\subset\Omega\setminus B_{\bar r}(0).$$
Then there exists $\bar x \in \Omega\setminus B_{\bar r}(0)$ such
that
$$u(\bar x )-kv(\bar x )-M=\max_{x\in\R^N\setminus\{0\}}(u-kv)(x)-M>0,$$
which follows by $\mathcal{A}_{k,M}\not=\O.$ Then, by Lemma \ref{cr
2.0}, we have
$$
(-\Delta)^\alpha (u-k v)(\bar x)\ge0,
$$
which is impossible with
\begin{eqnarray*}
(-\Delta)^\alpha (u-k v)(\bar x)&=&-u^p(\bar x )+kv^p(\bar x
)\\&\leq&-(k^p-k)v^p(x_0)-M^p
\\&<&0.
\end{eqnarray*}
We finish the proof.\hfill$\Box$

By the definition of $\mathcal{A}_{k,M}$ for any $M_1\ge M_2\ge0$,
we have that $\mathcal{A}_{k,M_1}\subset \mathcal{A}_{k,M_2}$.  For
notation convenient, we denote that
$\mathcal{A}_{k}=\mathcal{A}_{k,0}$.

\begin{lemma}\label{lm 3.2}Under the hypotheses of Theorem \ref{teo 1} part $(i)$,
 if
$$\mathcal{A}_{k}\not=\O,$$
where $k>1$ and $\mathcal{A}_{k}$ is given above. Then
 \begin{equation}\label{4.1.6}
 \lim_{r\to 0}\sup_{|x|=r}(u-kv)(x)=+\infty.
 \end{equation}
\end{lemma}
{\bf Proof.} If not, we have $\bar M:=\sup_{x\in
 \R^N\setminus\{0\}}(u-kv)(x)<+\infty.$ We see that $\bar M>0$ and there
doesn't exist point $\bar x $ achieving the supreme of $u-kv$ in $
\Omega\setminus\{0\}$. Indeed, if not, we can get a contradiction as
in the proof of Lemma \ref{lm 3.1}.

By Lemma \ref{lm 3.1}, $\mathcal{A}_{k}$ verifies (\ref{4.1.5}). Let
$x_0\in\mathcal{A}_{k}$ chosen later and $r=|x_0|/4$. In the
following, we will consider the function $$w_k=u-kv\ \ \rm{in}\ \
\R^N\setminus\{0\}.$$
 Under the hypotheses of Theorem \ref{teo 1} part $(i)$,
for all $x\in B_r(x_0)\cap\mathcal{A}_{k}$,
\begin{equation}\label{4.2.1}
 (-\Delta)^\alpha
 w_k(x)=-u^p(x)+kv^p(x),
\end{equation}
 then we have that
\begin{equation}\label{4.1.7} (-\Delta)^\alpha
w_k\le-K_1r^{\tau-2\alpha}\ \ \ \ \rm{in}\ \
B_r(x_0)\cap\mathcal{A}_{k}.
\end{equation}
where $\tau p=-\frac{2\alpha p}{p-1}=\tau-2\alpha$ and
$K_1=C^{-p}_0(k^{p}-k)>0$ with $C_0$ is from (\ref{4.1.3}).

We define $$w(x)=\frac{2\bar M}{V(0)} V(r(x-x_0)),\quad x\in \R^N,$$
 where $V$ is given in Lemma \ref{lm 2.2}, then we see that
\begin{equation}\label{4.1.8}
w(x_0)=\max_{x\in\R^N} w(x)=2\bar M
\end{equation}
and
\begin{equation}\label{4.1.9}
(-\Delta)^\alpha w\leq \frac{2\bar M}{V(0)} r^{-2\alpha}\ \ \ \
\rm{in}\ \ B_r(x_0).
\end{equation}
Let $x_0\in\mathcal{A}_{k}$  close enough to $0$ such that
$$\frac{2\bar M}{V(0)}\leq K_1r^\tau.$$
Combining (\ref{4.1.7}) with (\ref{4.1.9}), we have that
\begin{eqnarray*}
 (-\Delta)^\alpha (w_k+w)(x)\le0,\ \ \ x\in
 B_r(x_0)\cap\mathcal{A}_{k}.
\end{eqnarray*}
By Lemma \ref{cr 2.0} and $w_k(x_0)>0$, $w=0$ in $B^c_r(x_0)$, then
we have
\begin{equation}\label{ying}
\begin{array}{lll}
w(x_0)&<&w_k(x_0)+w(x_0)
\\[3mm]& \leq& \sup_{x\in B_r(x_0)\cap\mathcal{A}_{k}}(w_k+w)(x)
\\[3mm]&\le&\sup_{x\in (B_r(x_0)\cap\mathcal{A}_{k})^c}(w_k+w)(x)
\\[3mm]&=&\max\{\sup_{x\in B_r^c(x_0)}w_k(x),\sup_{x\in B_r(x_0)\cap \mathcal{A}_{k}^c}(w_k+w)(x)\}
\\[3mm]&\le&\max\{\bar M,\sup_{x\in B_r(x_0)\cap \mathcal{A}_{k}^c}(w_k+w)(x)\}.
\end{array}
\end{equation}
We first see the  contradiction in case of $\sup_{x\in B_r(x_0)\cap
\mathcal{A}_{k}^c}(w_k+w)>\bar M$. By $w_k\le 0$ in
$\mathcal{A}_{k}^c$, we have
$$\sup_{x\in B_r(x_0)\cap \mathcal{A}_{k}^c}(w_k+w)(x)\le\sup_{x\in
B_r(x_0)\cap \mathcal{A}_{k}^c}w(x),$$ which together with
(\ref{ying}), we have
$$w(x_0)<\sup_{x\in
B_r(x_0)\cap \mathcal{A}_{k}^c}w(x)\le \sup_{x\in
\R^N}w(x)=w(x_0),$$ which is impossible.

We finally see the  contradiction in case of $\sup_{x\in
B_r(x_0)\cap \mathcal{A}_{k}^c}(w_k+w)\le \bar M$. By (\ref{ying}),
we have
$$w(x_0)<\bar M.$$
which is impossible with (\ref{4.1.8}). We finish the proof. \hfill
$\Box$

\begin{remark}\label{re 3.02}
It is clear that
 $$\mathcal{A}_{k}\cap B_{d_0}\not=\O.$$
\end{remark}

\begin{remark}\label{re 3.2}
Let $m_k(t)=\max_{|x|=t}w_k(x)$, which is continuous in
$(0,+\infty)$ and $$m_k(t)=0,\quad\forall\ t\ge diam(\Omega).$$
Lemma \ref{lm 3.2} is equivalent to say: if there exist $t_0>0$ such
that
$$m_k(t_0)>0.$$ Then, $$\lim_{t\to0^+}m_k(t)=+\infty.$$
Moreover, let $m_0=\max_{t\in[t_0,\infty)}m_k(t)>0$, then for any
$C>m_0$, there exists $t_C\in(0,t_0)$ such that
$$m_k(t_C)=C\ \ \rm{and}\ \ \ m_k(t)\leq C\  \rm{for\ all}\  t\in[t_C,\infty). $$

\end{remark}

Directly by the results of Lemma \ref{lm 3.2}, we have
\begin{corollary}\label{co 3.1}
If $\mathcal{A}_{k}\not=\O$ with $k>1$, then
$\mathcal{A}_{k,M}\not=\O$ for any $M\ge0$.
\end{corollary}

\begin{lemma}\label{lm 3.3}
 Let $x_0\in
\mathcal{A}_{k}\cap B_{d_0}$, $r=|x_0|/4$ and $$Q_n=\{z\in  B_{\frac
rn}\ |\ w_k(z)>M_n\},\ \ n\in\N$$ with $M_n=\max_{\Omega\setminus
B_{\frac rn}(0)}w_k(x)$, then there exist $C_n>0$ $(n\ge1)$
independent of $x_0$ and $k$, such that
\begin{equation}\label{4.1.10}
 \lim_{n\to+\infty} C_n=0
\end{equation}
and
\begin{equation}\label{xuan1}
\int_{Q_n}\frac{w_k(z)-M_n}{|z-x|^{N+2\alpha}}dz\leq
C_nr^{\tau-2\alpha},\ \ \forall x\in B_r(x_0).
\end{equation}

\end{lemma}
{\bf Proof.}  By $v\ge0$ in $\R^N\setminus\{0\}$, $M_n\ge0$ and
(\ref{4.1.3}), we have
$$w_k(z)-M_n\le u(z)\le C_0|z|^\tau,\quad z\in B_{d_0}.$$
For $x\in B_r(x_0)$ with $r=|x_0|/4$ and $z\in Q_n$, we have
$$|x-z|\ge|x|-|z|\ge 3r-\frac rn>r.$$
Together with $Q_n\subset B_{\frac rn}\subset B_{d_0}$, we have
\begin{eqnarray*}
\int_{Q_n}\frac{w_k(z)-M_n}{|z-x|^{N+2\alpha}}dz
&\le&\int_{Q_n}\frac{u(z)}{|z-x|^{N+2\alpha}}dz
\\&\leq&C_0r^{-N-2\alpha}\int_{
 B_{\frac rn}}|z|^{\tau}dz\\
&\le&Cr^{-N-2\alpha}\int_0^{\frac rn}t^{\tau+N-1}dt\\
&\le&\frac C{n^{N+\tau}}r^{\tau-2\alpha}.
\end{eqnarray*}
Let $C_n=\frac C{n^{N+\tau}}$, then $\lim_{n\to+\infty}=0$.
The proof is complete. \hfill$\Box$\\

Now we give the proof of Theorem \ref{th 3.1} as follows:

\noindent{\bf Proof of Theorem \ref{th 3.1}.}  $ \mathcal{A}$ is
defined in (\ref{4.1.4}). If the conclusion of Theorem \ref{th 3.1}
under hypothesis $(i)$ in Theorem \ref{teo 1} isn't true, then
$\mathcal{A}\not=\O$.

 Let
$\bar x\in \mathcal{A}$ and $k_0\in (1,\frac{u(\bar x)}{v(\bar
x)})$. For example, $k_0=\frac{u(\bar x)+v(\bar x)}{2v(\bar x )}$.
We observe that $\bar x\in \mathcal{A}_{k_0}$.
 By Corollary \ref{co 3.1}, $\mathcal{A}_{k_0,1}$
 is open and nonempty. By Lemma \ref{lm 3.1}, we have that
\begin{equation}\label{M}
0\in\partial \mathcal{A}_{k_0,1}.
\end{equation}
By using Remark \ref{re 3.2}, there exists  $x_0\in
\mathcal{A}_{k_0,1}\cap B_{d_0}$ such that
$$u(x_0)-k_0 v(x_0)=\max_{x\in\Omega\setminus B_{4r}(0)}(u-k_0
v)(x),$$ where $r=|x_0|/4$. We recall that $w_{k_0}=u-k_0v$, then,
by (\ref{4.1.3}), for all $x\in B_r(x_0)\cap\mathcal{A}_{k_0,1}$, we
have
\begin{eqnarray*}
 (-\Delta)^\alpha
 w_{k_0}(x)&=&-u^p(x)+k_0v^p(x)
 \\&\le&-(k_0^{p}-k_0) v^p(x)
 \\&\le&-C^{-p}_0(k_0^{p}-k_0) |x|^{\tau p}
 \\&\le&-C^{-p}_0(k_0^{p}-k_0)(|x_0|-r)^{\tau p}
 \\&=:&-K_1r^{\tau-2\alpha},
\end{eqnarray*}
where $\tau=-\frac{2\alpha}{p-1}$,
$K_1=3^{\tau-2\alpha}C^{-p}_0(k_0^{p}-k_0)>0$ and $C_0$ is from
(\ref{4.1.3}). Then we have that
\begin{equation}\label{4.1.11} (-\Delta)^\alpha
w_{k_0}\le-K_1r^{\tau-2\alpha}\ \ \ \ \rm{in}\ \
B_r(x_0)\cap\mathcal{A}_{k_0,1},
\end{equation}
We redefine $$w(x)=\frac{K_1r^{\tau}}2 V(r(x-x_0))$$ for $x\in
\R^N,$ where $V$ is given in Lemma \ref{lm 2.2}, then we see that
\begin{equation}\label{4.1.12}
(-\Delta)^\alpha w\le\frac{K_1r^{\tau-2\alpha}}2\ \ \ \ \rm{in}\ \
B_r(x_0),
\end{equation}
Combining with (\ref{4.1.11}) and (\ref{4.1.12}), we have that
\begin{equation}\label{xuan}
 (-\Delta)^\alpha (w_{k_0}+w)(x)\le-\frac{K_1r^{\tau-2\alpha}}2,\ \ \ x\in
 B_r(x_0)\cap\mathcal{A}_{k_0,1}.
\end{equation}
Let
\begin{eqnarray*}
M_n:=\max_{x\in\overline{ B_{5r}\setminus B_{\frac rn}}}w_{k_0}(x),
\end{eqnarray*}
for $n\ge1$,  we have $x_0\in B_{5r}\setminus B_{\frac rn}$, then
 \begin{equation}\label{3.3.2}
 M_n\ge w_{k_0}(x_0)=
\max_{x\in\Omega\setminus B_{4r}(0)}(u-k_0 v)(x).
\end{equation}

 We denote that $$Q_n=\{z\in  B_{\frac rn}\ |\ w_{k_0}(z)>M_n\},\ n\in\N$$
and
\begin{equation} \label{3.2.5}
\bar w_n(x)=\left\{
\arraycolsep=1pt
\begin{array}{lll}
M_n,\ \ &\rm{if}\quad x\in Q_n,\\[2mm]
 (w_{k_0}+w)(x),\ \ \ \ &\rm{if\ not.}
\end{array}
\right.
\end{equation}
 By Lemma \ref{lm 3.3}, then there exists $n_0>1$ such
that
$$C_{n_0}\leq \frac{K_1}{2},$$
which, together with (\ref{xuan}), (\ref{xuan1}), we obtain
\begin{eqnarray*}
(-\Delta)^\alpha\bar w_{n_0}(x)&=&(-\Delta)^\alpha
(w_{k_0}+w)(x)+\int_{
Q_{n_0}}\frac{w_{k_0}(z)-M_{n_0}}{|z-x|^{N+2\alpha}}dz
\\&\le&-\frac{K_1}2r^{\tau-2\alpha}+C_{n_0}r^{\tau-2\alpha}
\\&\le &0,\quad \ x\in B_r(x_0)\cap\mathcal{A}_{k_0,1}.
\end{eqnarray*}
By Lemma \ref{cr 2.0}  and $w_{k_0}(x_0)>1$, $x_0\in
B_r(x_0)\cap\mathcal{A}_{k_0,1}$, $w=0$ in $B^c_r(x_0)$, then we
have
\begin{equation}\label{ying1}
\begin{array}{lll}
w(x_0)+1&<&w_{k_0}(x_0)+w(x_0)=\bar w_{n_0}(x_0)
\\[3mm]& \leq& \sup_{x\in B_r(x_0)\cap\mathcal{A}_{k_0,1}}\bar w_{n_0}(x)
\\[3mm]&\le&\sup_{x\in (B_r(x_0)\cap\mathcal{A}_{k_0,1})^c}\bar w_{n_0}(x)
\\[3mm]&=&\max\{\sup_{x\in B_r^c(x_0)}\bar w_{n_0}(x),\sup_{x\in B_r(x_0)\cap \mathcal{A}_{k_0,1}^c}\bar w_{n_0}(x)\}
\\[3mm]&\le&\max\{M_{n_0},\sup_{x\in B_r(x_0)\cap \mathcal{A}_{k_0,1}^c}(w_{k_0}+w)(x)\}.
\end{array}
\end{equation}
We first claim that  $\sup_{x\in B_r(x_0)\cap
\mathcal{A}_{k_0,1}^c}(w_{k_0}+w)\le M_{n_0}$. If not, by
$w_{k_0}\le 1$ in $\mathcal{A}_{k_0,1}^c$, we have
$$\sup_{x\in B_r(x_0)\cap \mathcal{A}_{k_0,1}^c}(w_{k_0}+w)(x)\le\sup_{x\in
B_r(x_0)\cap \mathcal{A}_{k_0,1}^c}w(x)+1,$$ which together with
(\ref{ying1}), we have
$$w(x_0)+1<\sup_{x\in
B_r(x_0)\cap \mathcal{A}_{k_0,1}^c}w(x)+1\le \sup_{x\in
\R^N}w(x)+1=w(x_0)+1,$$ which is impossible. Then we have that
\begin{equation}\label{ying2}
w_{k_0}(x_0)+w(x_0)\le\sup_{x\in B_r(x_0)\cap\mathcal{A}_{k_0,1}}
(w_{k_0}+w)(x)\le M_{n_0}.
\end{equation}

Since $\bar\Omega\setminus  B_{\frac r{n_0}}(0)$ is compact and
(\ref{3.3.2}), then there exists $x_1\in \overline{ B_{5r}\setminus
B_{\frac r{n_0}}}$ such that
$$w_{k_0}(x_1)=M_{n_0}.$$
Together with (\ref{ying2}), we have
\begin{equation}\label{ying3}
\frac{K_1V(0)}2r^\tau=w(x_0)< w_{k_0}(x_0)+w(x_0)\le M_{n_0}=
w_{k_0}(x_1).
\end{equation}
By (\ref{4.1.3}) and $\frac r{n_0}\le|x_1|\le 5r$, we have
$$\frac{K_1V(0)}2\frac{v(x_1)}{5^\tau C_0}\le \frac{K_1V(0)}2r^\tau\le
w_{k_0}(x_1)=u(x_1)-k_0v(x_1), $$ which implies that
\begin{equation}\label{4.2.4}
u(x_1)>(1+c_0)k_0v(x_1),
\end{equation} where
$$c_0=\frac{3^{\tau-2\alpha}(k_0^{p-1}-1) V(0)}{2C^2_0n_0^{-\tau}}>0.$$

 Now we repeat the  process above initiating by $x_1$. We know that $K_1$ is increasing with $k_0$,
 which is replaced by $k_1=(1+c_0)k_0$ and $n_0$ is independent of electing $x_0$,
 so we can keep our first choosing of $n_0$, then we have $x_2\in
 \mathcal{A}$ such that $$
u(x_2)>(1+c_1)k_1v(x_2)>(1+c_0)^2k_0v(x_2), $$ since
$$c_1=\frac{3^{\tau-2\alpha}(k_1^{p-1}-1) V(0)}{2C^2_0n_0^{-\tau}}>c_0.$$
 Proceeding  inductively,  we can
find a sequence $\{x_m\}\subset \mathcal{A}$ such that
$$u(x_m)>(1+c_0)^mk_0v(x_m),$$ which contradicts (\ref{4.1.3}).
 \hfill $\Box$\\


With the help of Theorem \ref{th 3.1}, we can prove Theorem \ref{teo
1}.

\noindent{\bf Proof the uniqueness in part $(i)$  of Theorem
\ref{teo 1}. } By $\mathcal{A}=\O$ in Theorem \ref{th 3.1}, then
$$u\le v\ \ \rm{in}\ \   B_{d_0}.$$
By using Theorem \ref{th 3.1} in domain $\{x\in B_{d_0}\ |\
u(x)<v(x)\}$, we see that
$$u\equiv v\ \ \rm{in}\ \   B_{d_0}.$$
Let $\tilde w:=u-v$ in $\R^N\setminus\{0\}$.

We first prove $\tilde w\ge0$ in $\R^N\setminus\{0\}$. If not, there
exists some point $\bar x \in\Omega\setminus B_{d_0}(0)$ such that
 $$\tilde w(\bar x)=\min_{x\in\R^N\setminus\{0\}}\tilde w(x)<0.$$ We observe, on the one hand, that
\begin{equation}\label{3.3}
(-\Delta)^\alpha \tilde w(\bar x )<0.
\end{equation}
  On the other hand, we have that
$$
 (-\Delta)^\alpha \tilde w(\bar x)=-u^{p}(\bar x )+v^{p}(\bar x )>0,
$$
which is impossible with (\ref{3.3}). By the same way, we get
$\tilde w\le0$ in $\R^N\setminus\{0\}$. Then we have that $u\equiv
v$ in $\R^N\setminus\{0\}$. We complete the proof.\hfill$\Box$


\setcounter{equation}{0}
\section{Nonexistence}
In this section, we focus on the nonexistence of classical solutions
under the hypotheses of Theorem \ref{teo 1} part $(iii)$. The idea
of the proof is as following: if there is a solution $u$ for
(\ref{eq1.1}) such that (\ref{1.70}) holds for some
$\tau\in(-N,0)\setminus\{2\alpha-N,-\frac{2\alpha}{p-1}\}$, there
exists some constants $C_2\ge C_1>0$ such that
$$C_1=\liminf_{x\to 0}u(x)|x|^{-\tau}\le\limsup_{x\to
0}u(x)|x|^{-\tau}= C_2.$$ We will find two sub solutions (or both
super solutions) $U_1$ and $U_2$ such that
$$\lim_{x\to 0}U_1(x)=\frac{C_1}2 ,\quad \lim_{x\to 0}U_2(x)=2C_2.$$ By
using Proposition \ref{th 4}  below, we will get a contradiction.
Therefore there is no solution under assumption of Theorem \ref{teo
1} part $(iii)$.

\begin{proposition}\label{th 4}
Under the hypotheses of Theorem \ref{teo 1} part $(iii)$, we suppose
that $U_1$ and $U_2$ are both  sub solutions (or both super
solutions) of (\ref{eq1.1})  and satisfy that $U_1=U_2=0$ in
$\Omega^c$ and
\begin{eqnarray*}
0&<&\liminf_{x\to 0}U_1(x)|x|^{-\tau}\le\limsup_{x\to
0}U_1(x)|x|^{-\tau}\\&<&\liminf_{x\to 0}U_2(x)|x|^{-\tau}\le
\limsup_{x\to 0}U_2(x)|x|^{-\tau}<+\infty,
\end{eqnarray*}
for some $\tau\in(-N,0)$. For the case $\tau p>\tau-2\alpha$, we
assume more that\\
$(i)$ in the case that $U_1,U_2$ are  sub solutions, there exist
$C>0$ and $\bar \delta>0$,
\begin{equation}\label{5.1}
(-\Delta)^\alpha U_2(x)\le -C|x|^{\tau-2\alpha},\quad x\in
B_{\bar\delta};
\end{equation}
or\\ $(ii)$ in the case that $U_1,U_2$ are  super solutions, there
exist $C>0$ and $\bar \delta>0$,
\begin{equation}\label{5.2}
(-\Delta)^\alpha U_1(x)\ge C|x|^{\tau-2\alpha},\quad x\in
B_{\bar\delta}.
\end{equation}
Then there doesn't exist any  solution $u$ of (\ref{eq1.1}) such
that
\begin{equation}\label{4.1}
\limsup_{x\to 0}\frac{U_1(x)}{u(x)}<1<\liminf_{x\to
0}\frac{U_2(x)}{u(x)}.
\end{equation}
\end{proposition}
\noindent{\bf Proof.} Here we only prove the case that $U_1$ and
$U_2$ are sub solutions of (\ref{eq1.1}) and the other case could be
done similarly.
 We prove it by contradiction.  Assume that there exists a solution $u$ for
(\ref{eq1.1}) satisfying (\ref{4.1}).  We observe that Lemma \ref{lm
3.1} and Lemma \ref{lm 3.3} hold in $\{x\in \Omega\ |\
u(x)-kU_1(x)>0\}$ for any $k>1$ and Lemma \ref{lm 3.2} holds for
$\{x\in \Omega\ |\ u(x)-kU_1(x)>1\}$.

Denote $\mathcal{C}_0=\{x\in \R^N\setminus\{0\}\ |\
U_2(x)>u(x)>U_1(x)>1\}$, which is open and nonempty by (\ref{4.1}).
By our hypothesis on $U_1$, $U_2$ and (\ref{4.1}), there exists
$C_0>1$ such that
\begin{equation}\label{4.3}
\frac1{C_0}\leq
U_1(x)|x|^{-\tau}<u(x)|x|^{-\tau}<U_2(x)|x|^{-\tau}\leq C_0,\ \ x\in
\mathcal{C}_0.
\end{equation}
Let  $\bar x \in \mathcal{C}_0$ and $k_0\in (1,\frac{U_2(\bar x
)}{u(\bar x )})$. We denote
$$\mathcal{C}_{k_0}:=\{x\in \R^N\setminus\{0\}\ |\ U_2(x)-k_0 u(x)>1\}$$
which, by Lemma \ref{lm 3.2}, is open and nonempty. By Lemma \ref{lm
3.1}, we have that
$$
0\in\partial \mathcal{C}_{k_0}.
$$
By using Remark \ref{re 3.2}, there exists  $x_0\in
\mathcal{C}_{k_0}$ such that
$$u(x_0)-k_0 v(x_0)=\max_{x\in\Omega\setminus B_{4r}(0)}(u-k_0
v)(x),$$ where $r=|x_0|/4$. Let $w_{k_0}=u-k_0U_1$. In the case of
$\tau p\le\tau-2\alpha$,   by (\ref{4.3}) we have
\begin{eqnarray*}
 (-\Delta)^\alpha
 w_{k_0}(x)&\leq&-u^p(x)+k_0U_1^p(x)
 \\&\le&-(k_0^{p}-k_0) U_1^p(x)
 \\&\le&-C^{-p}_0(k_0^{p}-k_0)(|x_0|-r)^{\tau p}
 \\&=:&-K_1r^{\tau-2\alpha},\quad x\in
 B_r(x_0)\cap\mathcal{C}_{k_0},
 \end{eqnarray*}
where $K_1=3^{\tau-2\alpha}C^{-p}_0(k_0^{p}-k_0)>0$ and $C_0$ is
from (\ref{4.3}).  In the case of $\tau p>\tau-2\alpha$,   by
(\ref{4.3}) and (\ref{5.1}) we have
\begin{eqnarray*}
 (-\Delta)^\alpha w_{k_0}(x)&\leq&-u^p(x)-Ck_0|x|^{\tau-2\alpha}
 \\&\le&-Ck_0|x|^{\tau-2\alpha},\quad x\in B_r(x_0)\cap\mathcal{C}_{k_0}.
 \end{eqnarray*}

 Proceeding as the  Proof of Theorem \ref{th 3.1}, we find a sequence $\{x_m\}\subset \mathcal{C}_0$ such that
$u(x_m)>(1+k_1)^mk_0U_1(x_m)$ for a certain constant $k_1>0$, which
contradicts (\ref{4.3}). Then there is no solution of (\ref{eq1.1})
satisfying (\ref{4.1}).
 \hfill $\Box$\\

Now we are in the position to prove  Theorem \ref{teo 1} part
$(iii)$.\\[1mm]
\noindent{\bf Proof of  Theorem \ref{teo 1} $(iii)$.} With the help
of Corollary \ref{coro 2.1},  for any given $t_1>t_2>0$, we
construct two sub solutions (or both super solutions) $U_1$ and
$U_2$ of (\ref{eq1.1}) such that
$$\lim_{x\to 0}U_1(x)|x|^{-\tau}=t_1,\quad \lim_{x\to 0}U_2(x)|x|^{-\tau}=t_2.$$
Then we use Proposition \ref{th 4}, we can get there is no solution
of (\ref{eq1.1}).

 We will prove
the nonexistence results in 3 cases.\\[1mm]
 \textbf{Case 1:}   $\tau\in(-N,-N+2\alpha)$
and  $\tau p>\tau-2\alpha$.  Denote that
$$W_{\mu,t}=tV_\tau-\mu \bar V\quad\rm{in}\ \ \R^N\setminus\{0\},$$
where $t,\mu>0$, $V_\tau$ is defined in (\ref{2.1}) and $\bar V$ is
the solution of (\ref{3.1}). By Corollary \ref{coro 2.1}$(i)$, for
$x\in B_{\delta_1}$, we have
\begin{eqnarray*}
 (-\Delta)^\alpha
 W_{\mu,t}(x)+
 |W_{\mu,t}|^{p-1}W_{\mu,t}(x)\leq-\frac{t}C|x|^{\tau-2\alpha}+t^p|x|^{\tau
 p}.
\end{eqnarray*}
 For any fixed $t>0$, there exists $\delta_2\in(0,\delta_1]$, for all
 $\mu\ge0$, we get
\begin{equation}\label{4.7}(-\Delta)^\alpha
 W_{\mu,t}(x)+
 |W_{\mu,t}|^{p-1}W_{\mu,t}(x)\leq0,\quad x\in B_{\delta_2}.\end{equation}
To consider $x\in \Omega\setminus B_{\delta_2}(0)$, in fact,
$(-\Delta)^\alpha V_\tau$ is bounded in $ \Omega\setminus
B_{\delta_2}(0)$ and
\begin{eqnarray*}
 (-\Delta)^\alpha
 W_{\mu,t}(x)+ |W_{\mu,t}|^{p-1}W_{\mu,t}(x)\leq C(t+t^p)-\mu,\quad x\in \Omega\setminus B_{\delta_2}(0).
\end{eqnarray*}
For given $t>0$, there exists $\mu(t)>0$  such that
\begin{equation}\label{4.9}
 (-\Delta)^\alpha
 W_{\mu(t),t}(x)+|W_{\mu,t}|^{p-1}W_{\mu,t}(x)\leq0,\ \ x\in\Omega\setminus B_{\delta_2}(0).
\end{equation}
Combining with (\ref{4.7}) and (\ref{4.9}), we have that for any
$t>0$, there exists $\mu(t)>0$ such that $$ (-\Delta)^\alpha
 W_{\mu(t),t}(x)+
|W_{\mu(t),t}|^{p-1}W_{\mu(t),t}(x)\leq0,\ \ \
x\in\Omega\setminus\{0\}.$$ For given $t_1>t_2>0$, there exist
$\mu(t_1)>0$ and $\mu(t_2)>0$ such that
\begin{eqnarray*}
t_2=\lim_{ x\to 0} W_{\mu(t_2),t_2}(x)|x|^{-\tau}<\lim_{x\to 0}
W_{\mu(t_1),t_1}(x)|x|^{-\tau}=t_1.
\end{eqnarray*}
Using Proposition \ref{th 4} with both sub solutions
$W_{\mu(t_1),t_1}$ and $ W_{\mu(t_2),t_2}$, there isn't any solution
$u$ of (\ref{eq1.1})
satisfying (\ref{1.70}).\\[1mm]
 \textbf{Case 2:}   $\tau\in(-N,-N+2\alpha)$
and $\tau p<\tau-2\alpha$.   We denote that
$$U_{\mu,t}=tV_\tau+\mu \bar V\quad \rm{in}\ \R^N\setminus\{0\},$$
where $t,\mu>0$. We know that $U_{\mu,t}>0$ in $\Omega$. By
Corollary \ref{coro 2.1} $(i)$, for $x\in B_{\delta_1}$,
\begin{eqnarray*}
 (-\Delta)^\alpha
 U_{\mu,t}(x)+
 U^p_{\mu,t}(x)\geq-Ct |x|^{\tau-2\alpha}+t^p|x|^{\tau p},
\end{eqnarray*}
for some $C>0$.
 For any fixed $t>0$, there exists $\delta_2\in(0,\delta_1]$, for all
 $\mu\ge0$, we have
\begin{equation}\label{4.70}(-\Delta)^\alpha
 U_{\mu,t}(x)+
 U^p_{\mu,t}(x)\geq0,\quad x\in B_{\delta_2}.\end{equation}
To consider $x\in \Omega\setminus B_{\delta_2}(0)$, in fact,
$(-\Delta)^\alpha V_\tau$ is bounded in $ \Omega\setminus
B_{\delta_2}(0)$ and
\begin{eqnarray*}
 (-\Delta)^\alpha
 U_{\mu,t}(x)+
 U^p_{\mu,t}(x)\geq-Ct+\mu,\quad x\in\Omega\setminus
B_{\delta_2}(0).
\end{eqnarray*}
For any given $t>0$, there exists $\mu(t)>0$  such that
\begin{equation}\label{4.80}
 (-\Delta)^\alpha
 U_{\mu(t),t}(x)+
 U^p_{\mu(t),t}(x)\geq0,\ \ x\in\Omega\setminus B_{\delta_2}(0).
\end{equation}
Combining with (\ref{4.70}) and (\ref{4.80}), we have that for any
$t>0$, there exists $\mu(t)>0$ such that $$ (-\Delta)^\alpha
 U_{\mu(t),t}(x)+
 U^p_{\mu(t),t}(x)\geq0,\ \ \ x\in\Omega\setminus\{0\}.$$
For given $t_1>t_2>0$, there exist $\mu(t_1)>0$ and $\mu(t_2)>0$
such that
\begin{eqnarray*}
t_2=\lim_{ x\to 0}U_{\mu(t_2),t_2}(x)|x|^{-\tau}<\lim_{ x\to
0}U_{\mu(t_1),t_1}(x)|x|^{-\tau}=t_1,
\end{eqnarray*}
Using Proposition \ref{th 4} with both super solutions
$U_{\mu(t_1),t_1}$ and $ U_{\mu(t_2),t_2}$, there isn't any solution of (\ref{eq1.1}) satisfying (\ref{1.70}).\\[1mm]
\textbf{Case 3:} $\tau\in(-N+2\alpha,0)$. By Corollary \ref{coro
2.1}$(ii)$, there exists
 $\delta_1>0$ such that
\begin{equation}\label{4.4}
(-\Delta)^{\alpha}V_\tau(x)>0,\ \ x\in B_{\delta_1}.
\end{equation}
Since $V_\tau$ is $C^2$  in $\Omega$, then there exists $C>0$ such
that
\begin{equation}\label{4.5}
|(-\Delta)^\alpha V_\tau(x)|\leq C,\ \ x\in\Omega\setminus
B_{\delta_1}(0).
\end{equation}
Let $\bar U:=V_\tau+C \bar V$, then we have $\bar U>0$ in $\Omega$
and
$$(-\Delta)^\alpha \bar U\ge 0\ \ \rm{in}\ \ \Omega.$$
Then, we have that $t\bar U$ is  super solution of (\ref{eq1.1}) for
any $t>0$. Using Proposition \ref{th 4}, there isn't any solution of
(\ref{eq1.1}) satisfying (\ref{1.70}). The proof is
complete.\hfill$\Box$\\

\noindent{\bf Proof of Remark \ref{re 2}.} Since $p\ge \frac
N{N-2\alpha}$, we have that $-\frac{2\alpha}{p-1}\le-N$. So by
Theorem \ref{teo 1}, it is only left to prove the case that
$\tau=\tau_0=2\alpha-N$. We denote that
$$U_{\mu,t}=tV_{\tau_0}+\mu \bar V\quad \rm{in}\ \R^N\setminus\{0\},$$
where $t,\mu>0$. We know that $U_{\mu,t}>0$ in $\Omega$. By
Corollary \ref{coro 2.1} $(iii)$, for $x\in B_{\delta_1}$,
\begin{eqnarray*}
 (-\Delta)^\alpha
 U_{\mu,t}(x)+
 U^p_{\mu,t}(x)\geq-Ct+t^p|x|^{\tau_0 p},
\end{eqnarray*}
for some $C>0$.
 For any fixed $t>0$, there exists $\delta_2\in(0,\delta_1]$, for all
 $\mu\ge0$, we have
\begin{equation}\label{4.70}(-\Delta)^\alpha
 U_{\mu,t}(x)+
 U^p_{\mu,t}(x)\geq0,\quad x\in B_{\delta_2}.
 \end{equation}
To consider $x\in \Omega\setminus B_{\delta_2}(0)$, in fact,
$(-\Delta)^\alpha V_\tau$ is bounded in $ \Omega\setminus
B_{\delta_2}(0)$ and
\begin{eqnarray*}
 (-\Delta)^\alpha
 U_{\mu,t}(x)+
 U^p_{\mu,t}(x)\geq-Ct+\mu,\quad x\in\Omega\setminus
B_{\delta_2}(0).
\end{eqnarray*}
For any given $t>0$, there exists $\mu(t)>0$  such that
\begin{equation}\label{4.80}
 (-\Delta)^\alpha
 U_{\mu(t),t}(x)+
 U^p_{\mu(t),t}(x)\geq0,\ \ x\in\Omega\setminus B_{\delta_2}(0).
\end{equation}
Combining with (\ref{4.70}) and (\ref{4.80}), we have that for any
$t>0$, there exists $\mu(t)>0$ such that $$ (-\Delta)^\alpha
 U_{\mu(t),t}(x)+
 U^p_{\mu(t),t}(x)\geq0,\ \ \ x\in\Omega\setminus\{0\}.$$
For given $t_1>t_2>0$, there exist $\mu(t_1)>0$ and $\mu(t_2)>0$
such that
\begin{eqnarray*}
t_2=\lim_{ x\to 0}U_{\mu(t_2),t_2}(x)|x|^{-\tau_0}<\lim_{ x\to
0}U_{\mu(t_1),t_1}(x)|x|^{-\tau_0}=t_1,
\end{eqnarray*}
Using Proposition \ref{th 4} with both super solutions
$U_{\mu(t_1),t_1}$ and $ U_{\mu(t_2),t_2}$, there isn't any solution
of (\ref{eq1.1}) satisfying (\ref{1.70}).
 The proof is
complete.\hfill$\Box$

\end{document}